\documentclass[a4paper,12pt]{amsart}

\usepackage{rotating,pdflscape,color,amssymb,hyperref,subfigure,psfrag,amsmath,eufrak,bbm,epsfig}
\usepackage{a4wide,amscd,hyperref}
\usepackage{tikz} 
 \usepackage[usenames,dvipsnames]{pstricks}
 \usepackage{pst-grad} 
 \usepackage{pst-plot} 
\usepackage[small]{caption}

\newcommand{\ld}{\ldots}
\newcommand{\beg}{\begin}
\newcommand{\en}{\end}

\newcommand{\trm}{\textrm}

\newcommand{\bgt}{\begin{itemize}}
\newcommand{\ent}{\end{itemize}}

\newcommand{\eqre}{\eqref}
\newcommand{\re}{\ref}
\newcommand{\la}{\label}

\newcommand{\si}{\sigma}

\newcommand{\p}{\mathbb{P}}

\newcommand{\Tr}{\operatorname{Tr}}

\newcommand{\Ninf}{\underset{N\to\infty}{\longrightarrow}}

\newcommand{\E}{\mathbb{E}}

\newcommand{\R}{\mathbb{R}}
\newcommand{\C}{\mathbb{C}}

\newcommand{\ud}{\mathrm{d}}

\newcommand{\pro}{probability }

\newcommand{\f}{\frac}
\newcommand{\ff}{\frac{1}}
\newcommand{\lf}{\left}
\newcommand{\ri}{\right}

\newcommand{\st}{such that }

\newcommand{\lam}{\lambda}

\newcommand{\ti}{\times}

\newcommand{\ste}{\, ;\, }
\newcommand{\mc}{\mathcal }

\newcommand{\eps}{\varepsilon}

\newcommand{\al}{\alpha}

\newcommand{\bbm}{\begin{bmatrix}}
\newcommand{\ebm}{\end{bmatrix}}
\newcommand{\bes}{\begin{equation*}}
\newcommand{\ees}{\end{equation*}}
\newcommand{\be}{\begin{equation}}
\newcommand{\ee}{\end{equation}}
\newcommand{\beqy}{\begin{eqnarray}}
\newcommand{\eeqy}{\end{eqnarray}}
\newcommand{\beq}{\begin{eqnarray*}}
\newcommand{\eeq}{\end{eqnarray*}}

\newcommand{\lto}{\longrightarrow}

\newcommand{\ie}{\emph{i.e. }}

\newcommand{\bpm}{\begin{pmatrix}}
\newcommand{\epm}{\end{pmatrix}}

\newcommand{\cd}{\cdots}

\newcommand{\bpr}{\beg{proof}}
\newcommand{\epr}{\en{proof}}

\newcommand{\del}{\delta}

\newcommand{\bX}{\mathbf{X}}
\newcommand{\bY}{\mathbf{Y}}

\newcommand{\ka}{\kappa}

\newcommand{\bi}{\mathbf{i}}

\newcommand{\bv}{\mathbf{v}}

\newcommand{\ba}{\mathbf{a}}
\newcommand{\bb}{\mathbf{b}}
\newcommand{\bz}{\mathbf{z}}

\long\def\symbolfootnote[#1]#2{\begingroup
\def\thefootnote{\fnsymbol{footnote}}\footnote[#1]{#2}\endgroup}



\author{Florent Benaych-Georges} \address{MAP 5, UMR CNRS 8145 - Universit\'e Paris Descartes, 45 rue des Saints-P\`eres 75270 Paris cedex~6, France.} \email{florent.benaych-georges@parisdescartes.fr}

\author{Sandrine P\'ech\'e} \address{LPMA, Universit\'e Paris Diderot, 5 rue Thomas Mann
75013 Paris.} 
\email{sandrine.peche@math.univ-paris-diderot.fr}

\title[]{Largest eigenvalues and  eigenvectors of band or sparse random matrices}

\begin{document}
\maketitle

\beg{abstract}In this text, we consider an random $N\ti N$ matrix $\bX$ \st all but $o(N)$ rows of $\bX$ have  $W$ non identically  zero entries, the other rows having lass than $W$ entries  (such as, for example, standard or cyclic band matrices). We always suppose that $1\ll W\ll N$.
We first prove that if the entries are independent, centered,  have variance  one, satisfy a certain tail upper-bound condition and $W\gg (\log N)^{6(1+\al)}$, where $\al$ is a positive parameter depending on the distribution of the entries, then the  largest eigenvalue of $\bX/\sqrt{W}$ converges to the upper bound of its limit spectral distribution, that is $2$, as for Wigner matrices. This   extends some previous results by Khorunzhiy and Sodin \cite{KHORUNZHIYSemProba,Sodin} where less hypotheses were made on $W$, but more hypotheses were made about the law of the entries and the structure of the matrix.  Then, under the same hypotheses, we prove a delocalization result for the eigenvectors of $\bX$. More precisely we show that eigenvectors associated to eigenvalues ``far enough" from zero  cannot be essentially localized on less than $W/\log(N)$ entries. This lower bound on the localization length has to be compared to the recent result by Steinerberger in \cite{Steinerberger}, which states that 
the localization length in the edge is $\ll W^{7/5}$ or 
there is strong interaction between two eigenvectors in an interval of
length $W^{7/5}$.\en{abstract}


\section{Introduction}
Random band matrices (\ie random Hermitian matrices with independent entries vanishing out of a band around the diagonal) have raised lots of attention recently. Indeed, varying the bandwidth $W$ from $1$ to the full size shows (in the large size limit)
a crossover between a strongly disordered  regime, with localized eigenfunctions and weak eigenvalue correlation,  and a weakly disordered regime, with extended eigenfunctions and strong eigenvalue repulsion.
It is conjectured (and explained  on a Physics level of rigor by Fyodorov and  Mirlin in \cite{FM}) that for Gaussian band matrices, 
 the \emph{localization strength} (\ie the typical number of coordinates bearing most of the $\ell^2$ mass) of a typical eigenvector in the bulk of the spectrum shall be of order $ {L\sim N\wedge W^2}$, so that eigenvectors of the bulk should be localized (resp. extended) if $W\ll \sqrt{N}$ (resp. $\gg \sqrt{N}$). The only rigorous result in the direction of localization is by Schenker  \cite{Schenker}. Therein it is proved that ${L\ll W^8}$ for Gaussian band matrices. On the other hand, delocalization in the bulk is proved by Erd\"os, Knowles,  Yau and  Yin \cite{EKYY} when $W\gg N^{4/5}$. In both regimes, it is known from Erd\"os and  Knowles \cite{EK2011CMP, ErdKnow} that typically (i.e. disregarding a negligible proportion of eigenvectors) ${L\gg W^{7/6}\wedge N}$ for a certain class of random band matrices (with sub-exponential tails and symmetric distribution).
We refer the reader to Spencer \cite{Spencer} and Erd\"os,  Schlein and  Yau \cite{ESY2} for a more detailed discussion on the localized/delocalized regime.

 Regarding the edges of the spectrum, little is known about the behavior of the extreme eigenvalues and the  typical localization length of the associated eigenvectors. 
 As far as the limit of the largest eigenvalue is concerned, Khorunzhiy proved in \cite{KHORUNZHIYSemProba}  that for Gaussian band matrices, if  $ (\log N)^{3/2}\ll W\ll N$, then the extreme eigenvalues converge to the bounds $u_{\pm}$ of the support of the limiting  spectral measure (which is the semicircle law). For matrices with cyclic band structure and Bernoulli entries,  Sodin extended this result to the case where $ \log N \ll W\ll N$  in  \cite{Sodin}, where he proved important results about the fluctuations of the extreme eigenvalues around their limits.  Concerning the localization length $L$ of the eigenvectors associated to the extreme eigenvalues, one could conjecture the following on the basis of the Thouless argument as explained in \cite{FM}. For eigenvectors associated to eigenvalues $\lam$ close to the bottom edge e.g. $u_{-}$, the localization strength should behave as  $L \sim N \wedge W^2(\lambda-u_{-})$.
Sodin's statement \cite{Sodin} combined with
Erd\"os-Knowles-Yau-Yin's results \cite{EKYY} suggest that this should hold true as soon as $W\gg N^{5/6}.$ 
Moreover, Steinerberger proved recently in \cite{Steinerberger} that for matrices with Bernoulli entries and cyclic band structure, with \pro tending to one, we have either $L\ll  W^{7/5}$ or there is strong interaction between two eigenvectors in an interval of
length $W^(7/5)$. Let us also mention that in the quite different framework of band matrices with heavy tailed entries,  a transition between the localized and the delocalized regime at the edge was proved by the authors of the present text in \cite{floSandrine}.

In this text, we consider an random $N\ti N$ Hermitian matrix $\bX$ \st rows of $\bX$ have $W$ non identically  zero entries  (such as, for example, standard or cyclic band matrices). We always suppose that $1\ll W\ll N$.
We first prove that if the entries are independent, centered,  have variance  one, satisfy a certain tail upper-bound condition and $W\gg (\log N)^{6(1+\al)}$, where $\al$ is a positive parameter depending on the distribution of the entries, then the  largest eigenvalue of $\bX/\sqrt{W}$ converges to the upper bound of its limit spectral distribution, that is $2$, as for Wigner matrices. This   extends the above mentioned   results by Khorunzhiy and Sodin \cite{KHORUNZHIYSemProba,Sodin} where less hypotheses were made on $W$, but more hypotheses were made about  the law of the entries (they use in a crucial way the fact that the entries are symmetrically distributed) and about the structure of the matrix (in our result, we only need that most rows have $W$ non zero entries, no matter what position the entries have on the matrix).  Then, under some close  hypotheses,   we prove a delocalization result for the eigenvectors of $\bX$, precisely that most of them cannot be essentially localized on less than $W/\log N $ entries. This lower bound on the localization length has to be compared to the recent result by Steinerberger in \cite{Steinerberger}, which states that 
the localization \emph{length} (here we use the word \emph{length} rather than \emph{strength} because in \cite{Steinerberger}, the author considers \emph{intervals} carrying most of the $\ell^2$-mass of the eigenvectors) in the edge is $\ll W^{7/5}$ or 
there is strong interaction between two eigenvectors in an interval of
length $W^{7/5}$.
The paper is organized as follows: our main results are stated in the next section, Theorem \re{CVnorme21012} is proved in Section \re{109131}, Theorem \re{anniDadAD1AD} is proved in Section \re{1091322}, and some technical results needed here are proved in Section \re{cc810121} and in the appendix.

\noindent{\bf Notation.} Here, $A\ll B$ means  that $A/B\lto 0$ as $N\to\infty$. Let $\|\bX\|$ denote the spectral radius of the Hermitian matrix $\bX$ and $\lam_{\max}(\bX)$ denote its largest eigenvalue.

 \section{Main results}

    We make the following hypotheses. 
     \beg{hyp}\la{hypentreeindep} The matrix $\bX=(X_{ij})$ is an $N\ti N$ Hermitian random matrix (depending implicitly on $N$) with independent entries (modulo the symmetry).   
     \en{hyp}
    \beg{hyp}\la{210121}There is $W=W(N)$ \st 
  \be\la{hypMinimaleW}1\ll W\ll N\ee and \st  on each row of $\bX$, the number of non identically zero entries is $\le W$,  with equality    on  all but $o(N)$ rows.
    All non identically zero entries of $\bX$ are centered with variance one. Moreover, there exist  constants $C\in [0,+\infty)$ and $\al\in [0,+\infty)$ \st for all $k\ge 2$, \be\la{hypMoments}\E[|X_{ij}|^k]\le (Ck)^{\al k},\ee uniformly on $N,i,j$.
    \en{hyp}

 Then the following theorem has been proved under weaker moment hypotheses in \cite{BMP91} (using the resolvent approach), but can easily be reproved here using a standard moment method as in \cite{alice-greg-ofer,bai-silver-book}.
 
 \beg{Th}\la{CVSCL}Under Hypotheses \re{hypentreeindep} and \re{210121}, as $N\to\infty$, the empirical spectral law of $\f{\bX}{\sqrt{W}}$ converges weakly in \pro towards the law $\ff{2\pi}\sqrt{4-x^2}\ud x$, with support $[-2,2]$.
 \en{Th}

 Our first result is the following one.  
 \beg{Th}\la{CVnorme21012}Under Hypotheses \re{hypentreeindep} and \re{210121},  suppose that $W$ satisfies\be\la{hypW} W\gg (\log N)^{6(1+\al)},\ee with $\al$ the constant of \eqre{hypMoments}. Then as $N\to\infty$, we have the convergence in \pro \be\la{2101221h2}\f{\lam_{\max}(\bX)}{\sqrt{W}}\lto 2.\ee
 \en{Th}
 \beg{rmk}\la{CVnorme21012rmk}{\rm  This theorem  extends some results of \cite{KHORUNZHIYSemProba} and \cite{Sodin}. In these papers, the convergence of \eqre{2101221h2} is proved under the respective hypotheses $W\gg (\log N)^{3/2}$ and $W\gg \log N$, but for some particular models of matrices: in \cite{KHORUNZHIYSemProba} the matrices considered are Gaussian and in \cite{Sodin}, they have Bernoulli distributed entries. Both make a crucial use of the fact that the entries are symmetrically distributed. Moreover, in both papers, the authors also suppose and rely heavily on a particular position of the non zero entries of $\bX$. We do not make such a hypothesis here.}
 \en{rmk}
 
 To state our main result, a lower bound on the  localization length of eigenvectors of $\bX$, we slightly modify the hypotheses.  \\
 Let $\bX$  be a random matrix satisfying Hypothesis \ref{hypentreeindep}. We make the following two assumptions.
    \beg{hyp}\la{3091212:46}  For a certain sequence $W=W(N)\gg1$, we have the convergence in \pro $\f{\|\bX\|}{2\sqrt{W}}\Ninf 1$. \en{hyp}
  For example, Hypothesis \re{3091212:46} is satisfied if $\bX$ satisfies the hypotheses of Theorem \re{CVnorme21012}  or those the papers \cite{KHORUNZHIYSemProba} and \cite{Sodin} (see  Remark \re{CVnorme21012rmk} above). 
However we emphasize that Hypothesis \ref{3091212:46}, focused on the extreme eigenvalues, does not make (at least directly) any assumption on the maximal number of non zero entries per row of $\bX$ (it may be $N$), neither on the relative growth of $W$ with respect to $N$.   
 
 We also reinforce the assumption on the tail of the distribution of the entries. 
 Let $C>0$ be fixed. 
   
  \beg{hyp}\la{3091212:462}The entries $X_{ij}$ belong to the set $\mc{E}_{C}$ of  complex  random variables $X$ \st $$\E X =0 ,\qquad \qquad \E  |X |^k \leq  (Ck)^{k/2},\: \forall k \in \mathbb{N}.$$\en{hyp}
Note that this   assumption reinforces  \eqref{hypMoments} as $\alpha\leq 1/2$. It is equivalent to the fact that there exists $\delta>0$ and $K>0$ such that 
 \be \label{gaussianint}\E e^{\del |X|^2}\le K.\ee With a slight abuse of notation, we denote $\mc{E}_{C}$ by $\mc{E}_{\del,K}$, as our proof mostly uses assumption (\ref{gaussianint}).
 
   

 The following theorem is the main result of this text. 
 
 We use the following Definition 7.1 from  Erd\"os, Schlein and Yau  \cite{ESY2}: 
for $L$   a positive integer and $\eta>0$, a unit vector $\bv=(v_1, \ld, v_N)\in \C^N$ is said to be   \emph{$(L, \eta)$-localized}  if    there exists a set $S\subset\{1,†\ld, N\}$  \st $|S|=L$ and $\sum_{j\in S^c}|v_j|^2\le \eta$. 

  \beg{rmk}{\rm The largest $L$ and $\eta$ are, the strongest the statement \emph{``there is no $(L,\eta)$-localized eigenvector"} is.}\en{rmk}
 
 \beg{Th}\la{anniDadAD1AD} We suppose Hypotheses \re{hypentreeindep},  \re{3091212:46} and \re{3091212:462}. Fix $\eta\in (0,1/2)$ and   choose $L=L(N)$ \st \be\la{3091212h5B} L\ll \f{W}{\log N}.\ee Let $\lam_1, \ld, \lam_N$ be the eigenvalues of $\bX$ and let $\bv_1, \ld, \bv_N$ be some associated normalized eigenvectors. Then for any $\ka$ \st $\sqrt{\eta/(1-\eta)}<\ka<1$,  $$\p(\exists i,\; |\lam_i|\ge 2\ka \sqrt{W}\trm{ and $\bv_i$ is $(L, \eta)$-localized})\;\Ninf \;0.$$
\en{Th}

      \beg{rmk}{\rm The same proof can also bring to a version of this theorem where $\eta=\eta(N)\lto 0$. In this case, $\ka=\ka(N)$ is allowed to tend to zero, thus the theorem allows to lower bound the localization length of most  eigenvectors of $\bX$.}\en{rmk}
    \beg{rmk}{\rm The estimate in Theorem  \ref{anniDadAD1AD} is almost sharp, as shown by the case where $\bX$ is the 
block diagonal matrix formed with $[N/W](+1)$ GUE matrices of size $W$ (or at most $W$ for the last block).}\en{rmk}
 
 \section{Proof of Theorem \re{CVnorme21012}}\la{109131}
 The proof goes along the same lines as the proof of  Theorem 2 in the paper \cite{KF} by F\"uredi and Koml\'os (see also Theorem 2.1.22 in \cite{alice-greg-ofer}). First note that by Theorem \re{CVSCL}, we already know that for any $\eta>0$, $\p(\lam_{\max}(\bX)<(2-\eta)\sqrt{W})\lto 0$.
 
 For any $\eta>0$, for any $k\ge 1$,  $$\p(\lam_{\max}(\bX)>(2+\eta)\sqrt{W})\le \p(\Tr \bX^{2k}\ge (2+\eta)^{2k}W^k)\le W^{-k}(2+\eta)^{-2k} \E\Tr \bX^{2k},$$ hence it suffices to find a sequence $k=k(N)$  \st for any $\eta>0$,  \be\la{2101221h} \E\Tr \bX^{2k} \ll W^{k}(2+\eta)^{2k}.\ee
 
We have $$\E\Tr \bX^{2k}=\sum \E X_{i_1i_2}\cd X_{i_{2k}i_{1}},$$ where the sum is over collections  $\bi=(i_1, \ld, i_{2k})$ \st for all $\ell$, $i_\ell\in \{1, \ld, N\}$.
For each $\bi$, let $G_{\bi}$ be the simple, non oriented graph with vertex set $\{i_1, \ld, i_{2k}\}$ and edges $\{i_{\ell}, i_{\ell+1}\}$ ($1\le \ell\le 2k$, with the convention $i_{2k+1}=i_1$). For
the expectation in the RHT above to be non zero, we need all edges to be visited at least twice by the path $\bi$ (because the $X_{ij}$'s are centered) and the edges $\{i_{\ell}, i_{\ell+1}\}$ to be \st $X_{i_\ell, i_{\ell+1}}$ is non identically zero. The symmetric group $S_N$ acts on the set of $\bi$'s by $\si\cdot(i_1, \ld, i_{2k}):=(\si(i_1), \ld, \si(i_{2k}))$. 
Following Section 2.1.3 of \cite{alice-greg-ofer}, we 
denote by $\mc{W}_{2k,t}$ the set of equivalence classes, under the action of $S_N$,  of $\bi$'s \st all edges of $G_{\bi}$ are   visited at least twice by the path $\bi$ and $G_{\bi}$ has exactly $t$ vertices (this set is actually stable under this action). 

Note that for $\mc{W}_{2k,t}$ to be non empty, we need to have $t\le k+1$. Indeed, $G_{\bi}$ is always connected hence its number of vertices is at most its number of edges plus one. 

Note that for any $w\in \mc{W}_{2k,t}$, the number of $\bi$'s in the class $w$ is at most $N W^{t-1}$. 

It follows from the previous remarks that  $$\E\Tr \bX^{2k}\le N \sum_{t=1}^{k+1}W^{t-1}\sum_{w\in\mc{W}_{2k,t}}\max_{\bi\in w}  \E X_{i_1i_2}\cd X_{i_{2k-1}i_{1}}.$$

Now, let us fix $t\in\{ 1, \ld, k+1\}$, $w\in  \mc{W}_{2k,t}$ and $\bi\in w$. Let us denote by $l$ (resp. $m$) the number of edges of $G_\bi$ visited exactly twice (resp. at least three times). Obviously, $2l+3m\le 2k$. 
Moreover, as $l+m$ is the number of edges of the $G_\bi$, hence by  connectedness of $G_\bi$ again, we have $t\le l+m+1$.   So 
$$6t\le 6m+6l+6=2(3m+2l)+2l+6\le 4k +2l+6,$$ so  $$2k-2l\le 6(k-t+1).$$

Now, notice that as the $X_{ij}$'s have variance one, $\E X_{i_1i_2}\cd X_{i_{2k}i_{1}}$ can be reduced to the expectation of a product of $2k-2l$ $X_{ij}$'s, hence by \eqre{hypMoments} and H\"older's inequality, $$\E X_{i_1i_2}\cd X_{i_{2k}i_{1}}\le \Big (C(2k-2l)\Big )^{\al (2k-2l)}\le \{6C(k-t+1)\}^{6\al (k-t+1)}.$$

As a consequence, $$\E\Tr \bX^{2k}\le N \sum_{t=1}^{k+1}W^{t-1}\#\mc{W}_{2k,t} \ti \{6C(k-t+1)\}^{6\al (k-t+1)}.$$
Now, we shall use Lemma 2.1.23 of \cite{alice-greg-ofer}, which states that
$\#\mc{W}_{2k,t}\le 4^k(2k)^{6(k-t+1)}$ as soon as $t\le k+1$ (the case $t=k+1$ is technically not contained in Lemma 2.1.23 of \cite{alice-greg-ofer}, but follows from Equation (2.1.20) and Lemma 2.1.3 of the same book).
It follows that \beq\E\Tr \bX^{2k}&\le& N  4^k\sum_{t=1}^{k+1}W^{t-1}(2k)^{6(k-t+1)}\{6C(k-t+1)\}^{6\al (k-t+1)}\\
&=& N W^{k}  4^k\sum_{i=0}^{k}W^{-i}(2k)^{6i}(6Ci)^{6\al i}\\
&\le & N W^{k}  4^k\lf(1-\f{(2k(6Ck)^\al)^6}{W}\ri)^{-1},
\eeq
where the last inequality is true as soon as $W>(2k(6Ck)^\al )^6$. Then it is easy to see that  \eqre{2101221h} holds for $k=k(N)$ \st $\log N\ll k\ll W^{\ff{6(1+\al)}}.$

      \section{Proof of Theorem \re{anniDadAD1AD}}\la{1091322}
   
   Before proving  Theorem \re{anniDadAD1AD},  we shall need  the following theorem and its corollary. The proof of Theorem  \re{2991219h9} is postponed to Section \re{cc810121}.
    
  \beg{Th}\la{2991219h9}Under  Hypotheses \re{hypentreeindep} and  \re{3091212:462}, there are constants $c_2=c_2(\del, K)>0$ and $C_2=C_2(\del, K)<\infty$ independent of all the other parameters \st for all $t>0$, \be\la{3091212:09}\p(\|\bX\|>t\sqrt{N})\le e^{-c_2(t^2-C_2)N}.\ee
  \en{Th}
  
  Let us denote by $\rho(\bX)$  the  spectral radius of $\bX$ and, for $L\ge 1$, by  $\rho_L(\bX)$   the maximum spectral radius of its $L\ti L$ principal submatrices (a \emph{principal} submatrix is a submatrix chosen by extracting a certain set of columns and \emph{the same} set of rows, but this set does not need to be an interval). 
 \beg{cor}\la{3091212h5A}Under  Hypotheses \re{hypentreeindep} and  \re{3091212:462},  there exists $t<\infty$ and $c_3>0$ \st \be\la{3091212h5}\p(\rho_L(\bX)\ge t\sqrt{L\log N})\le e^{-c_3L\log N}.\ee
 \en{cor}
 
 \bpr The number of ways to choose an $L\ti L$ principal submatrix is $\le N^L=e^{L\log N}$. 
 For each submatrix $\mathbf{S}$, $\p(\rho(\mathbf{S})\ge t\sqrt{N\log N})\le   \exp\{-c_2(t^2\log N-C_2)L\}$, hence by the union bound,  
 $$\p(\rho_L(\bX)\ge t \sqrt{L\log N})\le  \exp\{[-c_2(t^2\log N-C_2)+\log N]L\},$$ thus if $c_2t^2>1$, then \eqre{3091212h5} holds for a certain $c_3>0$. 
 \epr

To prove  Theorem \re{anniDadAD1AD},  we shall also  need  the following   lemma (see Lemma 4.2 in  \cite{floSandrine}).
 
 \beg{lem}For all $i$, if $\bv_i$ is $(L, \eta)$-localized, then 
$|\lam_i|\le \f{\rho_L(\bX)+\sqrt{\eta}\rho(\bX)}{\sqrt{1-\eta}}.$ 
  \en{lem}
  
  Let us now prove Theorem \re{anniDadAD1AD}.
 \bpr Let us choose $\eps>0$ \st  $\sqrt{\eta/(1-\eta)}(1+\eps)<\ka$ and set $$\del:=\ka-\sqrt{\eta/(1-\eta)}(1+\eps).$$ We know, by \eqre{3091212:46},  that with \pro tending to $1$, $\rho(\bX)\le (1+\eps)2\sqrt{W}$, \ie  $$\f{\sqrt{\eta}\rho(\bX)}{\sqrt{1-\eta}}\le (\ka-\del)2\sqrt{W}. $$ 
 Moreover, by Corollary \re{3091212h5A}, there is $t<\infty$ \st with \pro tending to one, $\rho_L(\bX)\le t\sqrt{L\log N}$. But by \eqre{3091212h5B},  for $N$ large enough, $$\f{t\sqrt{L\log N}}{\sqrt{1-\eta}}\le 2\del\sqrt{W},$$ so the theorem is proved.
 \epr

  \section{Proof of  Theorem  \re{2991219h9}}\la{cc810121}
  
  \subsection{A  preliminary lemma}

  
  \beg{lem}\la{lem309121}Under  Hypotheses \re{hypentreeindep} and  \re{3091212:462},  there are constants   $c_1, C_1$ depending only on $\del, K$ \st for any $\bz\in \C^N$ with $|\bz|\le 1$, \be\la{3091213h}\p(\bz^*\bX^*\bX\bz\ge Nt)\le e^{-c_1(t-C_1)N}.\ee \en{lem}
  
  \beg{rmk}{\rm If $\bX$ is a random $N\ti N$  matrix whose  maximum number of non identically zero entries per row is $W$ (like for a \emph{band matrix} with band width $W$), then \eqre{3091213h} remains true with $N$ replaced by $W$ everywhere (for some constants still depending only on $\del$ and $K$).}\en{rmk}
  
\bpr We denote by $\bX_1, \ld, \bX_N$ the columns of $\bX$.
 We have , for any $\tau,C$ as in Lemma \re{1791217:56} of the appendix, 
$$  \E e^{\tau^2 \bz^*\bX^*\bX\bz}\ =\  \E e^{\tau^2 |\bX\bz|^2}\
 =\ \E e^{\tau^2  \sum_j | \bX_j\cdot \bz|^2}\ \le \ e^{NC\tau^2}
$$
Hence 
   $$\p(\bz^*\bX^*\bX\bz\ge N t)\le  \E e^{\tau^2 \bz^*\bX^*\bX\bz}  e^{-\tau^2Nt}\le e^{-\tau^2(t-C)N}.$$ So the lemma is proved.\epr

    \subsection{Proof of Theorem \re{2991219h9}}
     
    \beg{lem}\la{igimyg2} Let $N\ge 1$. 
  For any fixed $0<\eps<1/4$, there exists a family $(\bz_i)_{i\in I}$ of elements of the unit ball of   $\C^N$ \st $|I|\le (2/\eps)^{2N}$ and any element of the unit sphere of $\C^N$ is within a distance at most $\eps$ of one of the $\bz_i$'s. Moreover, for any positive $N\ti N$ Hermitian matrix $P$, $$\lam_{\max}(P)\le \f{\max_i \bz_i^* P\bz_i}{1-2\eps}.$$ 
  \en{lem}
 The first part of this lemma is well known (see e.g. \cite{LedouxTalagrand}), whereas its second part  follows from the fact that for any vectors of the unit ball $\bz, \bz_i$, 
 $$\Big |  \bz_i^* P\bz_i- \bz^* P\bz \Big | \leq 2 ||P|| \ti ||\bz-\bz_i||$$ and specifying $\bz$ to be an eigenvector associated to $\lambda_{\max}.$
  Let us now prove Theorem \re{2991219h9}. 
  
  \bpr By Lemma \re{lem309121},  we know that there are constants there are constants $c_1, C_1$ depending only on $\del, K$ \st for any $\bz\in \C^N$ with $\|\bz\|\le 1$, $$\p(\bz^*\bX^*\bX\bz\ge Nt)\le e^{-c_1(t-C_1)N}.$$
   
     Now,  by Lemma \re{igimyg2}, we have \beq \p(\lam_{\max}(\bX^*\bX)\ge Nt)&\le& \p(\max_i \bz_i^* \bX^*\bX\bz_i\ge Nt(1-2\eps))\\ &\le& (2/\eps)^{2N} e^{-c_1(t(1-2\eps)-C_1)N}\\
  &=&  e^{(-c_1(1-2\eps)(t-\f{C_1}{1-2\eps})+2\log(2/\eps))N}\\
  &=&e^{-c_2(t-C_2)N}
  \eeq
  
 As a consequence, $$\p(\|\bX\|>t\sqrt{N})=\p(\lam_{\max}(\bX^*\bX)\ge t^2N)\le e^{-c_2(t^2-C_2)N}. $$
 \epr

 \section{Appendix: technical results}

  \beg{lem}\la{2891218:16}For any real centered random variable $X$ and any $r\in \R$, we have $$\E e^{rX}\le 1+3\E [e^{\del X^2}](e^{r^2/\del}-1)\le e^{3r^2\E [e^{ \del X^2}]/\del}$$ for any $\del>0$. 
  \en{lem}
  
  \bpr
   The second inequality follows from the fact that for any $y\ge 1$, we have the inequality  $1+y(e^{r^2/\del}-1)\le e^{r^2y/\del}$ (this is obvious with the series expansion of $\exp$).

    So let us prove the first inequality. 
  Note   that up to a replacement of $X$ by  $rX$ and of $\del$ by $\del/r^2$, we shall suppose that $r=1$. 
  
  The case where $\E[e^{\del X^2}]=\infty$ is obvious, hence we focus on the other case, which allows to expend all sums with the moments of $X$. 
  
  \emph{Claim :} for all $x\ge 0$, $e^x\le 1+x+3\f{e^x+e^{-x}-2}{2}$.
  Indeed, both terms are equal for $x=0$ and the derivative of $2$RHT$-2$LHT is $e^x-3e^{-x}+2,$ which is increasing, hence has the same sign as $x$. 
  
  It follows that $$\E e^X\le 1+3\sum_{n\ge 1}\f{\E X^{2n}}{(2n)!}\le 1+3\sum_{n\ge 1}\f{n!\del^{-n}\E e^{\del X^2}}{(2n)!}\le 1+3\E [e^{\del X^2}](e^{1/\del}-1),$$ where we used    first $\E X=0$, then    $\f{\del^{n}\E X^{2n}}{n!}\le \E e^{\del X^2}$ and at last $\f{n!}{(2n)!}\le \ff{n!}$.
  \epr
  
    Let $\ba\cdot\bb$ denote the standard scalar product of two complex vectors and let $|\cdot|$ denote the associated norm.
 \beg{lem}\la{1791217:56}  Let us fix $\del, K>0$. Then there is $\tau=\tau(\del, K)>0$ and $C=C(\del, K)>0$ \st for all $N\ge 1$, all $\bz\in \C^N$ \st $|\bz|\le \tau$, for any $\bY$   random vector taking values in $\C^N$ with independent  components in the set $\mc{E}_{\del,K}$ defined at Hypothesis \re{3091212:462}, $$\E e^{ |\bY\cdot \bz|^2}\le e^{C|\bz|^2}.$$
  \en{lem}
  
  \bpr \emph{First step:} Let us first prove the result for $\bY$ having independent components in  $$ \mc{E}_{\del,K}^{\R}:=\{Y\in  \mc{E}_{\del,K} \ste Y\trm{ is real-valued}\}$$ and $\bz\in \R^N$. 
  Let $\tau_\R>0$ be \st for any $t\in [0, \tau_\R)$, we have $$12t^2K/\del<1\qquad\trm{ and }\qquad \lf( 1-12t^2K/\del\ri)^{-1/2}\le e^{12t^2K/\del}.$$  
    Let $g$ be a standard real Gaussian variable, independent of the other variables, let $\E_g$ denote the expectation with respect to $g$ and let $\E$ denote the expectation with respect to all other variables than $g$. 

  For any $\tau>0$ and $\bz\in \R^N$ \st $|\bz|\le \tau$, using the formula $e^{x^2}=\E_g e^{\sqrt{2}xg}$, we have $$\E e^{|\bY\cdot \bz|^2} =\E \E_g e^{\sqrt{2} g\bY\cdot\bz }
  =\E_g \E e^{\sqrt{2} g\bY\cdot\bz }
  = \E_g\prod_{i} \E e^{\sqrt{2} g Y_iz_i}.$$ Hence by Lemma \re{2891218:16}, 
$$
\E e^{|\bY\cdot \bz|^2}  \le  \E_g\prod_{i} e^{6 g^2z_i^2K/\del}=  \E_g  e^{6 g^2|\bz|^2K/\del}=\lf( 1-12|\bz|^2K/\del\ri)^{-1/2} \le  e^{12|\bz|^2K/\del}
  .$$

 \emph{Second  step:}  To extend this result to the complex case, just decompose $Y$ and $\bz$ into real and imaginary parts and  use H\"older inequality to see that 
 the constants $\tau=\f{\tau_\R}{\sqrt{8}}$ and $C=4C_\R$ are convenient in the general case. 
  \epr

   \noindent{\bf Acknowledgments:} We would like to thank Stefan Steinerberger for having brought to our attention a misunderstanding we made from his paper \cite{Steinerberger} in the first version of this text. We would also like to thank the anonymous referee for some interesting remarks.

   \begin{thebibliography}{10}
   \bibitem{alice-greg-ofer}
G.~Anderson, A.~Guionnet, O.~Zeitouni
\newblock \emph{An Introduction to Random Matrices}.
\newblock  Cambridge studies in advanced mathematics, {118} (2009).
\bibitem{bai-silver-book} Z.D.~Bai, J.W.~Silverstein \emph{Spectral analysis of large dimensional random matrices}. Second Edition, Springer, New York, 2009.
  \bibitem{floSandrine} F. Benaych-Georges, S. P\'ech\'e \emph{Localization and delocalization  for heavy tailed band matrices}, to appear in Ann. Inst. Henri Poincar\'e Probab. Stat. 
  \bibitem{BMP91} L. V. Bogachev, S. A. Molchanov,  L. A. Pastur \emph{On the density of states of random
band matrices}, Mat. Zametki 50 (1991), 31--42, 157.
    \bibitem{EK2011CMP}  L. Erd\"os, A. Knowles
\emph{Quantum diffusion and eigenfunction delocalization in a random band matrix model}. 
Comm. Math. Phys. 303 (2011), no. 2, 509--554. 
 \bibitem{ErdKnow}  L. Erd\"os, A. Knowles \emph{ Quantum diffusion and delocalization for band matrices with general distribution.} Ann. Henri Poincaré\'e 12 (2011), no. 7, 1227--1319.   
  \bibitem{EKYY}  L. Erd\"os, A. Knowles, H-T. Yau, J. Yin \emph{Delocalization and Diffusion Profile for Random Band Matrices}.  arXiv:1205.5669, to appear in  {Comm. Math. Phys.}
 \bibitem{ESY2}L. Erd\"os, B. Schlein, H.T. Yau \emph{Semicircle law on short scales and delocalization of eigenvectors for Wigner random matrices}, Ann. Prob. 37 (2009), no. 3, 815--852.
 \bibitem{FM}Y.V. Fyodorov, A.D. Mirlin \emph{Scaling properties of localization in random band matrices: a $\sigma$-model approach}, Phys. Rev. Lett. 67 (1991), no. 18, 2405-2409.
 \bibitem{KF}Z. F\"uredi, J. Koml\'os  \emph{The eigenvalues of random symmetric matrices.}
Combinatorica 1 (1981), no. 3, 233-241.
 \bibitem{KHORUNZHIYSemProba} O. Khorunzhiy \emph{Estimates for moments of random matrices with Gaussian elements} in S\'eminaire de Probabilit\'es XLI, Lecture Notes in Math. 1934, Springer-Verlag, New York, 2008, pp. 51--92.
\bibitem{LedouxTalagrand} M. Ledoux, M. Talagrand \emph{Probability in Banach spaces}
 Classics in Mathematics. Springer (2011).
 \bibitem{Schenker} J. Schenker \emph{Eigenvector localization for random band matrices with power law band width.} Comm. Math. Phys. 290 (2009), no. 3, 1065-1097. 
  \bibitem{Sodin} S. Sodin \emph{The spectral edge of some random band matrices}, Ann. Math. 172 (2010), 2223--2251. 
  \bibitem{Spencer} T. Spencer \emph{Random banded and sparse matrices} The Oxford handbook of random matrix theory. Oxford University Press (2011) 471-488.
  \bibitem{Steinerberger} S. Steinerberger \emph{On Eigenvectors of Random Band Matrices with Large Band}, arXiv:1307.5753  \en{thebibliography}

\en{document}